\documentclass [10pt] {article}

 \usepackage [cp1251] {inputenc}
 \usepackage [english] {babel}
 \usepackage {indentfirst}
 \usepackage {amssymb}
 \usepackage {stmaryrd}
 \usepackage [mathscr] {eucal}
 \usepackage {mathtools}
 \usepackage {amsthm}
 \usepackage [curve] {xypic}

 \let \til\~         
 \let \gra\`
 \let \acu\'

 \def \Z{{\mathbb Z}}
 \def \Q{{\mathbb Q}}
 \def \R{{\mathbb R}}

 \def \E{{\boldsymbol{\mathcal E}}}            

 \def \F{{\mathscr F}} 

 \def \id{{\mathrm{id}}}
 \def \ins{{\mathrm{in}}}          
 \def \inc{{\mathrm{in}}}          
 \def \re{{\mathrm{rel}}}          

 \def \O{{\displaystyle\raisebox{4.2pt}{$\scriptscriptstyle\boldsymbol<$}\mkern-4.2mu{|}}}
 \def \o{{\scriptstyle\raisebox{2.3pt}{$\scriptscriptstyle<$}\mkern-4.5mu{|}}}
 \def \0{{\mathchoice\O\O\o\o}}

 \def \c{\times}       

 \def \<{\langle}      
 \def \>{\rangle}
 \def \`{\raisebox{1pt}{$\scriptscriptstyle<$}}            
 \def \'{\raisebox{1pt}{$\scriptscriptstyle>$}}

 \def \lsb{{[\![}}     
 \def \rsb{{]\!]}}
 \def \[{\lfloor}      
 \def \]{\rceil}
 \def \({\langle}                  
 \def \){\rangle}

 \DeclareMathOperator {\Hom} {Hom}

 \DeclareMathOperator {\rel} {rel} 

 \def \bou{\vee}       
 \def \Bou{\mathbin{\overline{\underline\vee}}}
 \def \cro{\times}     
 \def \sma{\wedge}     

 \def \le{\leqslant}
 \def \ge{\geqslant}
 \def \divs{\mathrel|} 

 \def \-{\overline}    
 \def \~{\widetilde}
 \def \^{\widehat}

 \def \+{\boldsymbol}  

 \def \xto{\xrightarrow}

 \renewcommand* {\%} [2]
 {\overset{#2}#1}

 \renewcommand* {\$} [3]
 {\overset{#2}{\underset{#3}#1}}

 \renewcommand* {\|} [2]
 {\mathbin{{#1}|_{#2}}}

 \newcommand* {\head} [1]
 {\subsubsection* {\mathversion{bold}#1}}

 \newcommand* {\subhead} [1]
 {\addvspace\medskipamount
  \noindent {\it #1\/}}

 \newenvironment* {claim} [1] []
 {\begin{trivlist}\item [\hskip\labelsep {\bf #1}] \it}
 {\end{trivlist} }

 \newenvironment* {cl} [1] []
 {\begin{trivlist}\item [\hskip\labelsep {\it #1}] \rm}
 {\end{trivlist} }

 \newenvironment* {demo} [1] []
 {\begin{trivlist}\item [\hskip\labelsep {\it #1}] }
 {\end{trivlist} }

 \newenvironment* {dem}
 {\begin{trivlist}\item}
 {\end{trivlist} }

 \makeatletter
 \def \QED{\displaymath@qed}
 \makeatother

\begin {document}

 \title {\large\bf
         Homotopy similarity of maps}

 \author {\normalsize\rm
          S.~S.~Podkorytov}

 \date {}

 \maketitle

 \vspace {-2\bigskipamount}
 
 \begin {abstract} \noindent
 Given based cellular spaces $X$ and $Y$,
 $X$ compact,
 we define a sequence of increasingly fine equivalences on
 the based-homotopy set $[X,Y]$.
 \end {abstract}


 \head {\S~1. Introduction}


 Let $X$ and $Y$ be based cellular spaces 
 ( = CW-complexes),
 $X$ compact.
 Let
 $Y^X$ be the set of based continuous maps
 $X\to Y$
 and
 $\<Y^X\>$ be the free abelian group associated with $Y^X$.
 An element $A\in\<Y^X\>$, an {\it ensemble}, has the form
 \begin {equation} \label {ens}
 A=
 \sum_i
 u_i\`a_i\',
 \end {equation}
 where
 $u_i\in\Z$ and
 $a_i\in Y^X$.
 Let $\F_r(X)$ be the set of {\it subspaces\/}
 ( = subsets containing the basepoint)
 $T\subseteq X$ containing at most $r$ points
 distinct from the basepoint.
 Introduce the subgroup
 $$
 \<Y^X\>^{(r+1)}=
 \{\,A:\text{$A|_T=0$ in $\<Y^T\>$ for all $T\in\F_r(X)$}\,\}\subseteq
 \<Y^X\>.
 $$
 We have
 $$
 \<Y^X\>=\<Y^X\>^{(0)}\supseteq\<Y^X\>^{(1)}\supseteq\dotso.
 $$
 For ensembles $A,B\in\<Y^X\>$,
 let
 $$
 A\%=rB
 $$
 mean that
 $B-A\in\<Y^X\>^{(r+1)}$.

 For maps $a,b\in Y^X$,
 we say that
 $a$ is {\it $r$-similar\/} to $b$,
 $$
 a\%\sim rb,
 $$
 when there exists an ensemble $A\in\<Y^X\>$
 given by \eqref{ens}
 with all $a_i\sim a$
 ($\sim$ denotes based homotopy)
 such that $A\%=r\`b\'$.
 A simple example is given in Section~3.

 Our main results state that
 the relation $\%\sim r$ is an equivalence
 (Theorem~8.1)
 and
 respects homotopy
 (Theorem~5.2).
 It follows that
 we get a sequence of
 increasingly fine equivalences on
 the based-homotopy set $[X,Y]$.

 We conjecture that,
 for $0$-connected $Y$,
 a map is $r$-similar to the constant map
 if and only if
 it lifts to the classifying space of
 the $(r+1)$th term of the lower central series of
 the loop group of $Y$.

 A related notion is that of
 a homotopy invariant of finite order
 \cite{sta,uns}.
 A function $f:[X,Y]\to L$,
 where $L$ is an abelian group,
 is called an invariant of {\it order\/} at most $r$
 when
 for any ensemble $A\in\<Y^X\>$
 given by \eqref{ens}
 the congruence $A\%=r0$ implies
 $$
 \sum_i
 u_if([a_i])=0.
 $$
 It is clear that
 $f([a])=f([b])$
 if
 $a\%\sim rb$
 and
 $f$ has order at most $r$.
 In \S~11,
 we give an example of two maps
 that
 are not $2$-similar
 but cannot be distinguished by invariants of order at most $2$.
 In the stable dimension range,
 invariants of order at most $r$
 were characterized
 in a way similar to
 our conjecture about $r$-similarity
 \cite{sta}.

 The relation between
 $r$-similarity and
 finite-order homotopy invariants
 is similar to that between
 $n$-equivalence and
 finite-degree invariants in knot theory
 \cite{Gus-1,Gus-2}.
 The example of \S~11
 is similar to
 that of \cite[Remark~10.8]{Gus-2}.


 \head {\S~2. Preliminaries}


 By a {\it space\/} we mean a based space
 (unless the contrary is stated explicitly).
 The basepoint of a cellular space is a vertex.
 The basepoint of a space $X$ is denoted by
 $\0_X$ or
 $\0$.
 A {\it subspace\/} contains the basepoint.
 A {\it cover\/} is a cover by subspaces.
 A {\it map\/} is a based continuous map.
 The constant map $X\to Y$ is denoted by
 $\0^X_Y$ or
 $\0$.
 A {\it homotopy\/} is a based homotopy.

 For a subspace $Z\subseteq X$,
 $\inc:Z\to X$ is the inclusion.
 A wedge of spaces comes with the insertions
 ( = coprojections):
 $$
 \ins_k:X_k\to X_1\bou\dotso\bou X_n.
 $$
 Maps $a_k:X_k\to Y$ form the map
 $$
 a_1\Bou\dotso\Bou a_n:
 X_1\bou\dotso\bou X_n\to
 Y.
 $$
 This notation is also used for homotopy classes.
 
 The formula $a\|\sim Zb$ means homotopy $a|_Z\sim b|_Z$.
 Similarly,
 equality of restrictions to a subset $C$ is denoted by $\|=C$.

 For a set $E$,
 the associated abelian group $\<E\>$ is freely generated by
 the elements $\`e\'$, $e\in E$.
 A function $t:E\to F$ between two sets induces
 the homomorphism
 $$
 \<t\>:\<E\>\to\<F\>,
 \qquad
 \`e\'\mapsto\`t(e)\'.
 $$

 For a cover $\Gamma$ of a space $X$,
 we put
 $$
 \Gamma(r)=
 \{\,\{\0\}\cup G_1\cup\dotso\cup G_s\subseteq X:
      G_1,\dotsc,G_s\in\Gamma,0\le s\le r\,\}.
 $$
 For ensembles $A,B\in\<Y^X\>$,
 the formula
 $$
 A\$=r\Gamma B
 $$
 means that
 $A\|=WB$ in $\<Y^W\>$
 for all $W\in\Gamma(r)$.

 Expressions with $?$ denote functions:
 for example,
 $?^2:\R\to\R$ is the function $x\mapsto x^2$.


 \head {\S~3. A simple example}


 Put $\E=\{0,1\}\subseteq\Z$.
 Fix $r\ge0$.
 For $d=(d_1,\dotsc,d_{r+1})\in\E^{r+1}$,
 put $|d|=d_1+\dotso+d_{r+1}$.
 Consider a wedge of spaces
 $$
 W=U_1\bou\dotso\bou U_{r+1}\bou V.
 $$
 Introduce the maps
 $$
 \Lambda(d)=
 \lambda_1(d_1)\bou\dotso\bou\lambda_{r+1}(d_{r+1})\bou\id_V:
 W\to W,
 \qquad
 d\in\E^{r+1},
 $$
 where the map $\lambda_k(e):U_k\to U_k$,
 for $e\in\E$,
 is
 $\id$
 if $e=1$ and
 $\0$
 if $e=0$.

 \begin {claim} [3.1. Lemma.]
 Let
 $X$ and $Y$ be spaces
 and
 $p:X\to W$ and $q:W\to Y$ be maps.
 Consider the ensemble $A\in\<Y^X\>$,
 $$
 A=\sum_{d\in\E^{r+1}}(-1)^{|d|}\`a(d)\',
 $$
 where
 $$
 a(d):
 X
 \xto{p}
 W
 \xto{\Lambda(d)}
 W
 \xto{q}
 Y.
 $$
 Then $A\%=r0$.
 \end {claim}

 \begin {demo} [Proof.]
 Take $T\in\F_r(X)$.
 There is a $k$ such that
 $p(T)\cap\ins_k(U_k)=\{\0_W\}$.
 Then $a(d)|_T$ does not depend on $d_k$.
 We get
 $$
 A|_T=
 \sum_{d\in\E^{r+1}}(-1)^{|d|}\`a(d)|_T\'=
 0.
 \QED
 $$
 \end {demo}

 \subhead {Example.}
 Consider the wedge
 $$
 W=S^{n_1}\bou\dotso\bou S^{n_{r+1}}
 $$
 ($n_1,\dotsc,n_{r+1}\ge1$).
 Put $m=n_1+\dotso+n_{r+1}-r$
 and
 let $p:S^m\to W$ be a map with
 $$
 [p]=\[\dotso\[[\ins_1],[\ins_2]\],\dotsc,[\ins_{r+1}]\]
 $$
 (the iterated Whitehead product)
 in $\pi_m(W)$.
 We show that
 $\0\%\sim rp$.
 Consider the maps
 $$
 a(d):
 S^m
 \xto{p}
 W
 \xto{\Lambda(d)}
 W,
 \qquad
 d\in\E^{r+1}.
 $$
 Put $1_{r+1}=(1,\dotsc,1)\in\E^{r+1}$.
 By Lemma~3.1,
 $$
 \sum_{d\in\E^{r+1}\setminus\{1_{r+1}\}}(-1)^{r-|d|}\`a(d)\'
 \%=r
 \`a(1_{r+1})\'.
 $$
 All $a(d)$ on the left side are homotopic to $\0$.
 On the right,
 $a(1_{r+1})=p$
 because $\Lambda(1_{r+1})=\id$.
 Thus
 $\0\%\sim rp$.


 \head {\S~4. Equipment of a cellular space}


 Let $Y$ be a compact {\it unbased\/} cellular space.
 In this section,
 we turn off our convention that
 {\it maps} and {\it homotopies} preserve basepoints.

 \begin {claim} [4.1. Lemma.]
 There exist homotopies
 $$
 q_t:Y^2\to Y
 \quad
 \text{and}
 \quad
 p_t:Y^2\to[0,1],
 \qquad
 t\in[0,1],
 $$
 such that
 \begin {equation} \label {bdy}
 q_0(z,y)=y,
 \qquad
 q_t(z,z)=z,
 \qquad
 p_0(z,y)=0,
 \qquad
 p_t(z,z)=t,
 \end {equation}
 and,
 for any $(z,y)\in Y^2$ and $t\in[0,1]$,
 one has
 \begin {equation} \label {alt}
 p_t(z,y)=0
 \quad
 \text{or}
 \quad
 q_t(z,y)=z.
 \end {equation}
 \end {claim}

 Roughly speaking,
 the inclusions $\{z\}\to Y$, $z\in Y$, form a parametric
 cofibration.

 \begin {demo} [Proof
                \rm (after {\cite[Exemple on p.~490]{Ser}}).]
 By \cite[Corollary~A.10]{Hat},
 $Y$ is an ENR.
 Embed it to $\R^n$
 and choose
 its neighbourhood $U\subseteq\R^n$
 and
 a retraction $r:U\to Y$.
 Choose $\epsilon>0$ such that
 $U$ includes all closed balls of radius $\epsilon$ with
 centres in $Y$.
 Consider the homotopy $l_t:(\R^n)^2\to\R^n$, $t\in[0,1]$,
 \begin {align*}
 l_t(z,y)&=y+\min(\epsilon t/|z-y|,1)(z-y),
 \qquad
 z\ne y, \\
 l_t(z,z)&=z.
 \end {align*}
 Put
 $$
 q_t(z,y)=r(l_t(z,y))
 \quad
 \text{and}
 \quad
 p_t(z,y)=\max(t-|z-y|/\epsilon,0).
 \QED
 $$
 \end {demo}

 \begin {claim} [4.2. Corollary.]
 One can continuously associate
 to each path $v:[0,1]\to Y$
 a homotopy $E_t(v):Y\to Y$, $t\in[0,1]$, such that
 $E_0(v)=\id$
 and
 $E_t(v)(v(0))=v(t)$.
 \end {claim}

 \begin {demo} [Proof.]
 Using Lemma~4.1,
 put
 $$
 E_t(v)(y)=
 \begin{cases}
 q_t(v(0),y) &
 \text{if $p_t(v(0),y)=0$}, \\
 v(p_t(v(0),y)) &
 \text{if $q_t(v(0),y)=v(0)$}.
 \end{cases}
 $$
 \vspace {-28pt}
 $$ \QED $$
 \end {demo}


 \head {\S~5. Coherent homotopies}


 Let $X$ and $Y$ be cellular spaces,
 $X$ compact.

 \begin {claim} [5.1. Lemma.]
 Consider
 an ensemble $A\in\<Y^X\>$,
 $$
 A=\sum_i
 u_i\`a_i\',
 $$
 and
 maps $b,\~b\in Y^X$,
 $b\sim\~b$.
 Then there exist maps $\~a_i\in Y^X$, $\~a_i\sim a_i$, such that
 the ensemble
 $$
 \~A=\sum_i
 u_i\`\~a_i\'
 $$
 has the following property:
 if $A\|=Z\`b\'$
 for some subspace $Z\subseteq X$,
 then $\~A\|=Z\`\~b\'$.
 \end {claim}

 \begin {demo} [Proof.]
 We have a homotopy $h_t\in Y^X$, $t\in[0,1]$, such that
 $h_0=b$
 and
 $h_1=\~b$.
 Replace $Y$ by a compact cellular subspace
 that includes the images of all $a_i$ and $h_t$.

 For $x\in X$,
 introduce the path $v_x=h_?(x):[0,1]\to Y$.
 We have
 $v_x(0)=b(x)$ and
 $v_x(1)=\~b(x)$.
 For a subspace $Z\subseteq X$,
 introduce the functions $e^Z_t:Y^Z\to Y^Z$, $t\in[0,1]$,
 $$
 e^Z_t(d)(x)=E_t(v_x)(d(x)),
 \qquad
 x\in Z,
 \quad
 d\in Y^Z,
 $$
 where $E_t$ is given by Corollary~4.2.
 For $d\in Y^Z$,
 we have the homotopy $e^Z_t(d)\in Y^Z$, $t\in[0,1]$.
 The diagram
 $$
 \xymatrix {
 Y^X
 \ar[rr]^-{e^X_t}
 \ar[d]_-{?|_Z} &&
 Y^X
 \ar[d]^-{?|_Z} \\
 Y^Z
 \ar[rr]^-{e^Z_t} &&
 Y^Z
 }
 $$
 is commutative.
 We have
 $e^Z_0=\id$
 because
 $$
 e^Z_0(d)(x)=
 E_0(v_x)(d(x))=
 d(x).
 $$
 We have
 $e^X_1(b)=\~b$
 because
 $$
 e^X_1(b)(x)=
 E_1(v_x)(b(x))=
 E_1(v_x)(v_x(0))=
 v_x(1)=
 \~b(x).
 $$
 Put $\~a_i=e^X_1(a_i)$.
 Since $a_i=e^X_0(a_i)$,
 we have $\~a_i\sim a_i$.
 We have
 $$
 (\`\~b\'-\~A)|_Z=
 \<e^X_1\>(\`b\'-A)|_Z=
 \<e^Z_1\>((\`b\'-A)|_Z).
 $$
 Thus
 $A\|=Z\`b\'$ implies $\~A\|=Z\`\~b\'$.
 \qed
 \end {demo}

 \begin {claim} [5.2. Theorem.]
 Let maps $a,b,\~a,\~b\in Y^X$ satisfy
 $$
 \~a\sim a\%\sim rb\sim\~b.
 $$
 Then $\~a\%\sim r\~b$.
 \end {claim}

 \begin {demo} [Proof.]
 By the definition of similarity,
 it suffices to show that
 $a\%\sim r\~b$.
 We have an ensemble $A\in\<Y^X\>$,
 $$
 A=\sum_i
 u_i\`a_i\',
 $$
 where $a_i\sim a$,
 such that
 $A\%=r\`b\'$.
 By Lemma~5.1,
 there is an ensemble $\~A\in\<Y^X\>$,
 $$
 \~A=\sum_i
 u_i\`\~a_i\',
 $$
 where $\~a_i\sim a_i$,
 such that
 $\~A\%=r\`\~b\'$.
 Since $a_i\sim a$,
 we have shown that
 $a\%\sim r\~b$.
 \qed
 \end {demo}


 \head {\S~6. Underlaying a cover}


 Let $X$ and $Y$ be cellular spaces,
 $X$ compact.

 \begin {claim} [6.1. Lemma.]
 Consider
 an ensemble $A\in\<Y^X\>$,
 $$
 A=\sum_i
 u_i\`a_i\'.
 $$
 Then there exist maps $\~a_i\in Y^X$, $\~a_i\sim a_i$, such that
 the ensemble
 $$
 \~A=\sum_i
 u_i\`\~a_i\'
 $$
 has the following property:
 if $A|_Z=0$
 for some subspace $Z\subseteq X$,
 then $\~A|_V=0$
 for some neighbourhood $V\subseteq X$ of $Z$.
 \end {claim}

 \begin {demo} [Proof.]
 Replace $Y$ by a compact cellular subspace
 that includes the images of all $a_i$.
 We will use the ``equipment'' $(q_t,p_t)$ given by Lemma~4.1.

 Let $i$ that numbers $a_i$ run over $1,\dotsc,n$.
 Define maps $a^k_i\in Y^X$, $1\le i\le n$, $0\le k\le n$,
 by the rules
 $a^0_i=a_i$
 and
 \begin {equation} \label {aki}
 a^k_i(x)=q_1(a^{k-1}_k(x),a^{k-1}_i(x)),
 \qquad
 x\in X,
 \end {equation}
 for $k\ge1$.
 Put $\~a_i=a^n_i$.
 We have $a^k_i\sim a^{k-1}_i$ because
 $a^k_i=h_1$ and
 $a^{k-1}_i=h_0$
 for the homotopy $h_t\in Y^X$, $t\in[0,1]$,
 $$
 h_t(x)=q_t(a^{k-1}_k(x),a^{k-1}_i(x)),
 \qquad
 x\in X.
 $$
 Thus $\~a_i\sim a_i$.

 \begin {cl} [Claim 1.]
 If $a^{k-1}_i\|=Qa^{k-1}_j$
 for some subspace $Q\subseteq X$,
 then $a^k_i\|=Qa^k_j$.
 \end {cl}

 \begin {dem}
 This follows from \eqref{aki}.
 \end {dem}
 
 \begin {cl} [Claim 2.]
 If $a^{i-1}_i\|=Qa^{i-1}_j$
 for some subspace $Q\subseteq X$,
 then there exists
 a neighbourhood $W\subseteq X$ of $Q$ such that
 $a^i_i\|=Wa^i_j$.
 \end {cl}

 \begin {dem}
 Indeed,
 if $a^{i-1}_i\|=Qa^{i-1}_j$,
 then,
 by~\eqref{bdy},
 $$
 p_1(a^{i-1}_i(x),a^{i-1}_j(x))=1
 $$
 for $x\in Q$.
 There exists a neighbourhood $W\subseteq X$ of $Q$ such that
 $$
 p_1(a^{i-1}_i(x),a^{i-1}_j(x))>0
 $$
 for $x\in W$.
 Then,
 by~\eqref{alt},
 $$
 q_1(a^{i-1}_i(x),a^{i-1}_j(x))=a^{i-1}_i(x)
 $$
 for $x\in W$.
 By \eqref{aki},
 $$
 a^i_i(x)=
 q_1(a^{i-1}_i(x),a^{i-1}_i(x))=
 a^{i-1}_i(x)
 $$
 (because $q_1(z,z)=z$ by \eqref{bdy})
 and
 $$
 a^i_j(x)=
 q_1(a^{i-1}_i(x),a^{i-1}_j(x)).
 $$
 Thus
 $a^i_i(x)=a^i_j(x)$
 for $x\in W$,
 as required.
 \end {dem}

 Take a subspace $Z\subseteq X$.

 \begin {cl} [Claim 3.]
 If $a_i\|=Za_j$,
 then there exists
 a neighbourhood $W\subseteq X$ of $Z$ such that
 $\~a_i\|=W\~a_j$.
 \end {cl}

 \begin {dem}
 This follows from
 the construction of $\~a_i$
 and
 the claims 1 and 2.
 \end {dem}

 Consider the equivalence
 $$
 R=\{\,(i,j):a_i\|=Za_j\,\}
 $$
 on the set $I=\{1,\dotsc,n\}$.
 It follows from the claim~3 that
 there exists a neighbourhood $V\subseteq X$ of $Z$ such that
 $\~a_i\|=V\~a_j$
 for all $(i,j)\in R$.
 We have the commutative diagram
 $$
 \xymatrix {
 Y^Z &&
 I
 \ar[ll]_-{a_i|_Z\mapsfrom i:l}
 \ar[rr]^-{d:i\mapsto\~a_i|_V}
 \ar[d]^-{\pi} &&
 Y^V \\
 &&
 I/R,
 \ar[llu]^-{\-l}
 \ar[rru]_-{\-d} &&
 }
 $$
 where
 $\pi$ is the projection.
 The function $\-l$ is injective.
 Consider the elements
 $U\in\<I\>$,
 $$
 U=\sum_iu_i\`i\',
 $$
 and
 $\-U=\<\pi\>(U)\in\<I/R\>$.
 We have
 $$
 A|_Z=
 \<l\>(U)=
 \<\-l\>(\-U)
 \quad
 \text{and}
 \quad
 \~A|_V=
 \<d\>(U)=
 \<\-d\>(\-U).
 $$
 If $A|_Z=0$,
 then $\-U=0$
 because $\<\-l\>$ is injective.
 Then
 $\~A|_V=0$.
 \qed
 \end {demo}

 \begin {claim} [6.2. Corollary.]
 Consider
 an ensemble $A\in\<Y^X\>$,
 $$
 A=\sum_i
 u_i\`a_i\',
 $$
 such that
 $A\%=r0$.
 Then there exist maps $\~a_i\in Y^X$,
 $\~a_i\sim a_i$,
 such that
 the ensemble
 \begin {equation} \label {tildeA-0}
 \~A=\sum_i
 u_i\`\~a_i\'
 \end {equation}
 satisfies the condition
 $\~A\$=r\Gamma0$
 for some open cover $\Gamma$ of $X$.
 \end {claim}

 \begin {demo} [Proof.]
 Since $A\%=r0$,
 we have $A\|=T0$
 for all $T\in\F_r(X)$.
 By Lemma~6.1,
 there are maps $\~a_i\in Y^X$,
 $\~a_i\sim a_i$,
 such that the ensemble $\~A$ given by \eqref{tildeA-0}
 satisfies the condition $\~A\|={V(T)}0$
 for some neighbourhood $V(T)\subseteq X$ of
 each $T\in\F_r(X)$.
 There is an open cover $\Gamma$ of $X$ such that
 every $W\in\Gamma(r)$ is included in $V(T)$
 for some $T\in\F_r(X)$.
 Then
 $\~A\|=W0$
 for all $W\in\Gamma(r)$,
 that is,
 $\~A\$=r\Gamma0$.
 \qed
 \end {demo}

 \begin {claim} [6.3. Lemma.]
 Consider
 an ensemble $A\in\<Y^X\>$,
 $$
 A=\sum_i
 u_i\`a_i\',
 $$
 and
 a map $b\in Y^X$.
 Then there exist maps $\~a_i\in Y^X$,
 $\~a_i\sim a_i$,
 such that
 the ensemble
 \begin {equation} \label {tildeA}
 \~A=\sum_i
 u_i\`\~a_i\'
 \end {equation}
 has the following property:
 if $A\|=Z\`b\'$
 for some subspace $Z\subseteq X$,
 then $\~A\|=V\`b\'$
 for some neighbourhood $V\subseteq X$ of $Z$.
 \end {claim}

 \begin {demo} [Proof.]
 Let $\Pi$ be the set of subspaces $Z\subseteq X$ such that
 $A\|=Z\`b\'$.
 By Lemma~6.1,
 there are maps $\-a_i,\-b\in Y^X$,
 $\-a_i\sim a_i$ and
 $\-b\sim b$,
 such that
 the ensemble
 $$
 \-A=\sum_i
 u_i\`\-a_i\'
 $$
 satisfies the condition $\-A\|={V(Z)}\`\-b\'$
 for some neighbourhood $V(Z)\subseteq X$ of each $Z\in\Pi$.
 By Lemma~5.1,
 there are maps $\~a_i\in Y^X$,
 $\~a_i\sim\-a_i$,
 such that
 the ensemble $\~A$ given by \eqref{tildeA}
 satisfies the condition $\~A\|={V(Z)}\`b\'$
 for all $Z\in\Pi$.
 \qed
 \end {demo}

 \begin {claim} [6.4. Corollary.]
 Consider
 an ensemble $A\in\<Y^X\>$,
 $$
 A=\sum_i
 u_i\`a_i\',
 $$
 and
 a map $b\in Y^X$.
 Suppose that
 $A\%=r\`b\'$.
 Then there exist maps $\~a_i\in Y^X$,
 $\~a_i\sim a_i$,
 such that
 the ensemble
 \begin {equation} \label {tildeA-1}
 \~A=\sum_i
 u_i\`\~a_i\'
 \end {equation}
 satisfies the condition
 $\~A\$=r\Gamma\`b\'$
 for some open cover $\Gamma$ of $X$.
 \end {claim}

 \begin {dem}
 This follows from Lemma~6.3 as
 Corollary~6.2 does from Lemma~6.1.
 \qed
 \end {dem}

 \head {\S~7. Symmetric characterization of similarity}


 Let $X$ and $Y$ be cellular spaces,
 $X$ compact.

 \begin {claim} [7.1. Lemma.]
 Consider
 a cover $\Gamma$ of $X$,
 an open subspace $G\in\Gamma$,
 a closed subspace $D\subseteq X$,
 $D\subseteq G$,
 and
 maps $a,b_0,b_1\in Y^X$ such that
 $a\|\sim Gb_0$,
 $b_0\sim b_1\rel X\setminus D$,
 and
 $a\$\sim{r-1}\Gamma b_0$ in the following sense:
 there is an ensemble $A\in\<Y^X\>$,
 $$
 A=\sum_i
 u_i\`a_i\',
 $$
 where $a_i\sim a$,
 such that
 $A\$={r-1}\Gamma\`b_0\'$.
 Then
 there exists an ensemble $C\in\<Y^X\>$,
 $$
 C=\sum_k
 w_k\`c_k\',
 $$
 where $c_k\sim a$,
 such that
 $C\%=r\`b_1\'-\`b_0\'$.
 \end {claim}

 \begin {demo} [Proof.]
 There is a homotopy $h_t\in Y^X$, $t\in[0,1]$, such that
 $h_s=b_s$, $s=0,1$,
 and
 $h_t\|={X\setminus D}b_0$.
 Choose
 a continuous function $\phi:X\to[0,1]$ such that
 $\phi|_E=1$
 and
 $\phi|_{X\setminus F}=0$
 for some subspaces $E,F\subseteq X$,
 $E$ open,
 $F$ closed,
 such that
 $$
 D\subseteq
 E\subseteq
 F\subseteq
 G.
 $$

 Let $p\in Y^G$ be a map such that
 $p\sim b_0|_G$.
 Choose a homotopy $K_t(p)\in Y^G$, $t\in[0,1]$, such that
 $K_0(p)=p$,
 $K_1(p)=b_0|_G$,
 and, moreover,
 $K_t(p)=b_0|_G$
 if $p=b_0|_G$.
 Define a homotopy $L_t(p)\in Y^G$, $t\in[-1,1]$, by the rules
 $$
 L_t(p)(x)=
 K_{\phi(x)(t+1)}(p)(x),
 \qquad
 x\in G,
 $$
 for $t\in[-1,0]$
 and
 $$
 L_t(p)(x)=
 \begin {cases}
 h_t(x) &
 \text{if $x\in E$}, \\
 K_{\phi(x)}(p)(x) &
 \text{if $x\in G\setminus D$}
 \end {cases}
 $$
 for $t\in[0,1]$.
 We have
 $L_{-1}(p)=p$,
 $L_s(p)\|=Eb_s$, $s=0,1$,
 $L_0(p)\|={G\setminus D}L_1(p)$,
 and
 $L_t(p)\|={G\setminus F}p$.
 Moreover,
 $L_s(b_0|_G)=b_s|_G$,
 $s=0,1$.

 Let $d\in Y^X$ be a map such that
 $d\|\sim G b_0$.
 Define a homotopy $l_t(d)\in Y^X$, $t\in[-1,1]$,
 by the rules
 $l_t(d)\|=G L_t(d|_G)$
 and
 $l_t(d)\|={X\setminus F}d$.
 We have
 $l_{-1}(d)=d$,
 $l_s(d)\|=Eb_s$, $s=0,1$,
 $l_0(d)\|={X\setminus D}l_1(d)$,
 and
 $l_t(d)\|={X\setminus F}d$.

 Since $a_i\sim a\|\sim Gb_0$,
 the homotopies $l_t(a_i)$ are defined.
 Put
 $$
 C=
 \sum_i
 u_i(\`l_1(a_i)\'-\`l_0(a_i)\').
 $$
 We have $l_s(a_i)\sim a_i\sim a$.
 It remains to show that
 $C\%=r\`b_1\'-\`b_0\'$.
 Take $T\in\F_r(X)$.
 We check that
 \begin {equation} \label {tocheck}
 C
 \|=T
 \`b_1\'-\`b_0\'.
 \end {equation}
 We are in one of the following three cases.

 \begin {demo} [Case 1: \rm $T\cap D=\{\0_X\}$.]
 We have
 $l_0(a_i)\|=Tl_1(a_i)$
 and
 $b_0\|=Tb_1$.
 Thus
 both the sides of \eqref{tocheck} are zero on $T$.
 \end {demo}

 \begin {demo} [Case 2: \rm $T\cap F=\{\0_X,x_*\}$,
                where $x_*\in E$ and $x_*\ne\0_X$.]
 Put $Z=T\setminus\{x_*\}$.
 We have
 $Z\in\F_{r-1}(X)$ and
 $Z\cap F=\{\0_X\}$.
 Define functions $e_s:Y^Z\to Y^T$, $s=0,1$,
 by the rules
 $e_s(q)|_Z=q$ and
 $e_s(q)(x_*)=b_s(x_*)$.
 We have
 $e_s(b_0|_Z)=b_s|_T$ and
 $e_s(a_i|_Z)=l_s(a_i)|_T$.
 Thus
 $$
 \big(
 \`b_0\'-
 \sum_i
 u_i\`a_i\'
 \big)
 \big|_Z
 \xmapsto{\<e_s\>}
 \big(
 \`b_s\'-
 \sum_i
 u_i\`l_s(a_i)\'
 \big)
 \big|_T.
 $$
 Since $A\%={r-1}\`b_0\'$,
 the expression on the left is zero.
 Thus
 the one on the right is also zero,
 which implies \eqref{tocheck}.
 \end {demo}

 For a finite space $Z$,
 let $\Vert Z\Vert$ be
 the cardinality of $Z\setminus\{\0\}$.

 \begin {demo} [Case 3: \rm $\Vert T\cap G\Vert\ge2$.]
 We have $T=W\cup Z$
 for some subspaces $W,Z\subseteq X$ such that
 $W\cap Z=\{\0_X\}$,
 $W\subseteq G$,
 $Z\cap F=\{\0_X\}$,
 and
 $\Vert Z\Vert\le r-2$.
 Consider the subspace $M=G\cup Z\subseteq X$.
 Define functions $f_s:Y^M\to Y^T$, $s=0,1$.
 Take $q\in Y^M$.
 If $q\|\sim Gb_0$,
 put
 $f_s(q)\|=WL_s(q|_G)$ and
 $f_s(q)\|=Zq$.
 Otherwise,
 put $f_s(q)=\0^T_Y$.
 We have
 $f_s(b_0|_M)=b_s|T$
 and
 $f_s(a_i|_M)=l_s(a_i)|_T$.
 Thus
 $$
 \big(
 \`b_0\'-
 \sum_i
 u_i\`a_i\'
 \big)
 \big|_M
 \xmapsto{\<f_s\>}
 \big(
 \`b_s\'-
 \sum_i
 u_i\`l_s(a_i)\'
 \big)
 \big|_T.
 $$
 Since
 $M$ is included in some element of $\Gamma(r-1)$
 and
 $A\$={r-1}\Gamma\`b_0\'$,
 the expression on the left is zero.
 Thus
 the one on the right is also zero,
 which implies \eqref{tocheck}.
 \qed
 \end {demo}
 \end {demo}

 \begin {claim} [7.2. Lemma.]
 Let $a,b,\~b\in Y^X$ be maps such that
 $a\%\sim{r-1}b\sim\~b$
 and
 (*)
 $a\|\sim Sb$
 for any $S\in\F_1(X)$.
 Then
 there exists an ensemble $C\in\<Y^X\>$,
 $$
 C=\sum_k
 w_k\`c_k\',
 $$
 where $c_k\sim a$,
 such that
 $C\%=r\`\~b\'-\`b\'$.
 \end {claim}

 The condition (*) is satisfied automatically
 if $X$ or $Y$ is $0$-connected.
 It also follows from the condition $a\%\sim{r-1}b$
 if $r\ge2$
 (cf.\ the proof of Theorem~7.3).

 \begin {demo} [Proof.]
 There is an ensemble $A\in\<Y^X\>$,
 $$
 A=\sum_i
 u_i\`a_i\',
 $$
 whre $a_i\sim a$,
 such that
 $A\%={r-1}\`b\'$.
 Using Corollary~6.4,
 replace each $a_i$ by a homotopic map
 to get $A\$={r-1}\Gamma\`b\'$
 for some open cover $\Gamma$ of $X$.

 Call a subspace $G\subseteq X$ {\it primitive\/} if
 the map $\inc:G\to X$ is homotopic to the composition
 $$
 G\xto{f}
 S\xto{\inc}
 X
 $$
 for some
 subspace $S\in\F_1(X)$
 and
 map $f$.
 Since
 $X$ is Hausdorff and locally contractible,
 for any
 open subspace $U\subseteq X$
 and
 point $x\in U$,
 there exists a primitive open subspace $G\subseteq X$ such that
 $x\in G$
 and
 $G\subseteq U$.
 We replace the cover $\Gamma$ by its refinement
 consisting of primitive open subspaces.
 Then
 it follows from (*) that
 $a\|\sim G b$
 for each $G\in\Gamma$.

 Choose a finite partition of unity subordinate to $\Gamma$:
 $$
 \sum_{j=1}^m
 \phi_j
 =
 1,
 $$
 where
 each $\phi_j:X\to[0,1]$ is a continuous function such that
 $\phi_j|_{X\setminus D_j}=0$
 for some closed subspace $D_j\subseteq X$ such that
 $D_j\subseteq G_j$
 for some $G_j\in\Gamma$.
 Choose a homotopy $h_t\in Y^X$, $t\in[0,1]$, such that
 $h_0=b$
 and
 $h_1=\~b$.
 Define maps $b_j\in Y^X$, $0\le j\le m$, by the rule
 $$
 b_j(x)=h_{\phi_1(x)+\dotso+\phi_j(x)}(x).
 $$
 We have
 $b_0=b$,
 $b_m=\~b$,
 and
 $b_{j-1}\sim b_j\rel X\setminus D_j$.

 Take $j\ge1$.
 Applying Lemma~5.1 to
 the congruence $A\$={r-1}\Gamma\`b\'$ and
 the homotopy $b\sim b_{j-1}$,
 we get an ensemble $A_j\in\<Y^X\>$,
 $$
 A_j=\sum_i
 u_i\`a_{ji}\',
 $$
 where $a_{ji}\sim a_i$ ($\sim a$),
 such that
 $A_j\$={r-1}\Gamma\`b_{j-1}\'$.
 We have $a\|\sim{G_j}b\sim b_{j-1}$.
 By Lemma~7.1,
 there is an ensemble $C_j\in\<Y^X\>$,
 $$
 C_j=
 \sum_k
 w_{jk}\`c_{jk}\',
 $$
 where $c_{jk}\sim a$,
 such that
 $C_j\%=r\`b_j\'-\`b_{j-1}\'$.

 We get
 $$
 \sum_{j=1}^m
 C_j=
 \`b_m\'-\`b_0\'=
 \`\~b\'-\`b\'.
 \QED
 $$
 \end {demo}

 \begin {claim} [7.3. Theorem.]
 Consider
 maps $a,b\in Y^X$
 and
 ensembles $A,B\in\<Y^X\>$,
 $$
 A=
 \sum_i
 u_i\`a_i\'
 \quad
 \text{and}
 \quad
 B=
 \sum_j
 v_j\`b_j\',
 $$
 where
 $$
 \sum_i
 u_i=
 \sum_j
 v_j=
 1,
 $$
 $a_i\sim a$, and
 $b_j\sim b$,
 such that
 $A\%=rB$.
 Then
 $a\%\sim rb$.
 \end {claim}

 \begin {demo} [Proof.]
 Induction on $r$.
 If $r\le0$,
 the assertion is trivial.
 Suppose $r\ge1$.

 For $S\in\F_1(X)$,
 we have $a\|\sim Sb$
 because
 $$
 \`[a|_S]\'=
 \sum_i
 u_i\`[a_i|_S]\'=
 \lsb A|_S\rsb=
 \lsb B|_S\rsb=
 \sum_j
 v_j\`[b_j|_S]\'=
 \`[b|_S]\'
 $$
 in $\<[S,Y]\>$.
 Here
 $\lsb?\rsb:\<Y^S\>\to\<[S,Y]\>$ is the homomorphism induced by
 the projection $[?]:Y^S\to[S,Y]$.

 By induction hypothesis,
 $a\%\sim{r-1}b$.
 Take $j$.
 Since $b\sim b_j$,
 Lemma~7.2 gives an ensemble $C_j\in\<Y^X\>$,
 $$
 C_j=\sum_k
 w_{jk}\`c_{jk}\',
 $$
 where $c_{jk}\sim a$,
 such that
 $C_j\%=r\`b_j\'-\`b\'$.
 We have
 $$
 A-
 \sum_j
 v_jC_j\%=r
 A-
 \sum_j
 v_j(\`b_j\'-\`b\')=
 A-B+\`b\'\%=r
 \`b\',
 $$
 which proves the assertion.
 \qed
 \end {demo}


 \head {\S~8. Similarity is an equivalence}


 Let $X$ and $Y$ be cellular spaces,
 $X$ compact.

 \begin {claim} [8.1. Theorem.]
 The relation $\%\sim r$ on $Y^X$ is an equivalence.
 \end {claim}

 This was conjectured by A.~V.~Malyutin.

 \begin {demo} [Proof.]
 Reflexivity is trivial.
 Symmetry follows from Theorem~7.3.
 It remains to prove transitivity.

 Let maps $a,b,c\in Y^X$ satisfy $a\%\sim rb\%\sim rc$.
 There are ensembles $A,B\in\<Y^X\>$,
 $$
 A=
 \sum_i
 u_i\`a_i\'
 \quad
 \text{and}
 \quad
 B=
 \sum_j
 v_j\`b_j\',
 $$
 where
 $a_i\sim a$
 and
 $b_j\sim b$,
 such that
 $A\%=r\`b\'$ and
 $B\%=r\`c\'$.
 For each $j$,
 we have $b\sim b_j$ and,
 by Lemma~5.1,
 there is an ensemble $A_j\in\<Y^X\>$,
 $$
 A_j=
 \sum_i
 u_i\`a_{ji}\',
 $$
 where $a_{ji}\sim a_i$ ($\sim a$),
 such that
 $A_j\%=r\`b_j\'$.
 We have
 $$
 \sum_j
 v_jA_j\%=r
 \sum_j
 v_j\`b_j\'=
 B\%=r
 \`c\'.
 $$
 Thus
 $a\%\sim rc$.
 \qed
 \end {demo}

 Using Theorem~5.2,
 we introduce the relation of $r$-similarity on $[X,Y]$:
 $$
 [a]\%\sim r[b]
 \quad
 \Leftrightarrow
 \quad
 a\%\sim rb.
 $$
 It follows from Theorem~8.1 that
 it is an equivalence.


 \head {\S~9. The Hopf invariant}


 Let $X$ and $Y$ be spaces.
 Let
 $e\in C^m(Y)$ and
 $f\in C^n(Y)$
 ($m,n\ge1$)
 be (singular) cocycles and
 $g\in C^{m+n-1}(Y)$ be a cochain with $\delta g=ef$.
 Put
 $$
 [X,Y]_{e,f}=
 \{\,\+a:\text{$\+a^*([e])=0$ and
               $\+a^*([f])=0$ in $H^\bullet(X)$}\,\}
 \subseteq
 [X,Y]
 $$
 and
 $$
 Y^X_{e,f}=
 \{\,a:[a]\in[X,Y]_{e,f}\,\}
 \subseteq
 Y^X.
 $$
 Given $a\in Y^X_{e,f}$,
 choose a cochain $p\in C^{m-1}(X)$ such that
 $\delta p=a^\#(e)$
 and
 put
 $$
 q=
 pa^\#(f)-a^\#(g)\in
 C^{m+n-1}(X).
 $$
 Then
 $\delta q=0$
 and
 the class $[q]\in H^{m+n-1}(X)$
 neither depends on the choice of $p$
 nor changes if
 $a$ is replaced by a homotopic map.
 Putting $h([a])=[q]$,
 we get the function
 $$
 h:[X,Y]_{e,f}\to H^{m+n-1}(X),
 $$
 which we call the {\it Hopf invariant\/}
 \cite{Ste}.

 \begin {claim} [9.1. Lemma.]
 Let
 $X_0$ be a space
 and
 $t:X\to X_0$
 be a map.
 We have the Hopf invariants
 $$
 h_0:[X_0,Y]_{e,f}\to H^{m+n-1}(X_0)
 \quad
 \text{and}
 \quad
 h:[X,Y]_{e,f}\to H^{m+n-1}(X).
 $$
 Given $a_0\in Y^{X_0}$,
 put $a=a_0\circ t\in Y^X$.
 If $a_0\in Y^{X_0}_{e,f}$,
 then
 $a\in Y^X_{e,f}$
 and
 $h([a])=t^*(h_0([a_0]))$ in $H^{m+n-1}(X)$.
 \qed
 \end {claim}

 \begin {claim} [9.2. Lemma.]
 Take elements
 $\+u\in\pi_m(Y)$ and
 $\+v\in\pi_n(Y)$.
 Put
 $$
 \Delta=
 \(\+u^*([e]),[S^m]\)
 \(\+v^*([f]),[S^n]\)+
 (-1)^{mn}
 \(\+u^*([f]),[S^m]\)
 \(\+v^*([e]),[S^n]\)
 \in\Z
 $$
 (the last two Kronecker indices vanish
 unless $m=n$).
 Consider
 the Hopf invariant
 $$
 h:
 [S^{m+n-1},Y]_{e,f}\to
 H^{m+n-1}(S^{m+n-1})
 $$
 and the Whitehead product $\[\+u,\+v\]\in\pi_{m+n-1}(Y)=[S^{m+n-1},Y]$.
 Then
 $\[\+u,\+v\]\in[S^{m+n-1},Y]_{e,f}$ and
 $$
 \(h(\[\+u,\+v\]),[S^{m+n-1}]\)=
 (-1)^{mn+m+n}\Delta.
 $$
 \end {claim}

 Caution:
 the sign in the last equality is sensitive to certain conventions.

 \begin {demo} [Proof \rm (after {\cite[\S~19]{Ste}}).]
 We assume that
 $S^m\bou S^n\subseteq S^m\cro S^n$ in the standard way.
 We have the commutative diagram
 $$
 \xymatrix {
 S^{m+n-1}
 \ar[r]^-{\phi}
 \ar[d]_-{\inc} &
 S^m\bou S^n
 \ar[d]^-{\inc} \\
 D^{m+n}
 \ar[r]^-{\chi} &
 S^m\cro S^n,
 }
 $$
 where $[\phi]=\[[\ins_1],[\ins_2]\]$ in $\pi_{m+n-1}(S^m\bou S^n)$.
 We have the chain of
 homomorphisms and
 sendings
 \begin {equation} \label {homol}
 \xymatrix {
 H_{m+n-1}(S^{m+n-1}) &
 \scriptstyle{[S^{m+n-1}]} \\
 H_{m+n}(D^{m+n},S^{m+n-1})
 \ar[u]^-{\partial}
 \ar[d]_-{(\chi,\phi)_*} &
 \scriptstyle{[D^{m+n}]}
 \ar@{|->}[u]
 \ar@{|->}[d] \\
 H_{m+n}(S^m\cro S^n,S^m\bou S^n) &
 \scriptstyle{\re_*([S^m\cro S^n])} \\
 H_{m+n}(S^m\cro S^n).
 \ar[u]^-{\re_*} &
 \scriptstyle{[S^m\cro S^n]}
 \ar@{|->}[u]
 }
 \end {equation}

 Choose representatives
 $u:S^m\to Y$ and
 $v:S^n\to Y$
 of
 $\+u$ and
 $\+v$,
 respectively.
 Consider the maps
 $$
 a:
 S^{m+n-1}
 \xto{\phi}
 S^m\bou S^n
 \xto{w=u\Bou v}
 Y.
 $$
 Clearly,
 $[a]=\[\+u,\+v\]$
 in $\pi_{m+n-1}(Y)$.

 Choose
 cocycles
 $\^e\in C^m(S^m\cro S^n)$ and
 $\^f\in C^n(S^m\cro S^n)$ and
 a cochain $\^g\in C^{m+n-1}(S^m\cro S^n)$
 such that
 $$
 \^e|_{S^m\bou S^n}=w^\#(e),
 \qquad
 \^f|_{S^m\bou S^n}=w^\#(f),
 \quad
 \text{and}
 \quad
 \^g|_{S^m\bou S^n}=w^\#(g).
 $$

 We have
 $$
 a^\#(e)=
 \phi^\#(w^\#(e))=
 \phi^\#(\^e|_{S^m\bou S^n})=
 \chi^\#(\^e)|_{S^{m+n-1}}
 $$
 in $C^m(S^{m+n-1})$.
 It follows that
 $a^*([e])=0$
 in $H^m(S^{m+n-1})$
 (which is automatic
 unless $n=1$).
 Similarly,
 $a^*([f])=0$
 in $H^n(S^{m+n-1})$.
 Thus
 $[a]\in[S^{m+n-1},Y]_{e,f}$.

 Let $z_k\in H^k(S^k)$ be the class with $\(z_k,[S^k]\)=1$.
 One easily sees that
 $$
 [\^e]=
 \(\+u^*([e]),[S^m]\)(z_m\times1)+
 \(\+v^*([e]),[S^n]\)(1\times z_n)
 $$
 in $H^m(S^m\cro S^n)$ and
 $$
 [\^f]=
 \(\+v^*([f]),[S^n]\)(1\times z_n)+
 \(\+u^*([f]),[S^m]\)(z_m\times1)
 $$
 in $H^n(S^m\cro S^n)$.
 Thus
 $[\^e][\^f]=\Delta(z_m\times z_n)$
 in $H^{m+n}(S^m\cro S^n)$ and
 \begin {equation} \label {int}
 \([\^e][\^f],[S^m\cro S^n]\)=
 (-1)^{mn}\Delta.
 \end {equation}

 Choose a cochain $\~p\in C^{m-1}(D^{m+n})$ such that
 $\delta\~p=\chi^\#(\^e)$.
 Put
 $$
 \~q=\~p\chi^\#(\^f)-\chi^\#(\^g)
 \in C^{m+n-1}(D^{m+n}).
 $$
 Put
 $$
 p=\~p|_{S^{m+n-1}}
 \in C^{m-1}(S^{m+n-1})
 \quad
 \text{and}
 \quad
 q=\~q|_{S^{m+n-1}}
 \in C^{m+n-1}(S^{m+n-1}).
 $$
 We have
 $$
 \delta p=
 \delta\~p|_{S^{m+n-1}}=
 \chi^\#(\^e)|_{S^{m+n-1}}=
 \phi^\#(\^e|_{S^m\bou S^n})=
 \phi^\#(w^\#(e))=
 a^\#(e)
 $$
 and
 \begin {multline*}
 q=
 p\chi^\#(\^f)|_{S^{m+n-1}}-\chi^\#(\^g)|_{S^{m+n-1}}=
 p\phi^\#(\^f|_{S^m\bou S^n})-\phi^\#(\^g|_{S^m\bou S^n})= \\ =
 p\phi^\#(w^\#(f))-\phi^\#(w^\#(g))=
 pa^\#(f)-a^\#(g).
 \end {multline*}
 Thus
 $\delta q=0$
 and
 $h([a])=[q]$.

 We have
 $$
 \delta\~q=
 \chi^\#(\^e)\chi^\#(\^f)-\delta\chi^\#(\^g)=
 \chi^\#(\^e\^f-\delta\^g).
 $$
 We have the chain of
 homomorphisms and
 sendings
 $$
 \xymatrix {
 H^{m+n-1}(S^{m+n-1})
 \ar[d]_-{\delta} &
 \scriptstyle{[q]}
 \ar@{|->}[d] \\
 H^{m+n}(D^{m+n},S^{m+n-1}) &
 \scriptstyle{[\chi^\#(\^e\^f-\delta\^g)]} \\
 H^{m+n}(S^m\cro S^n,S^m\bou S^n)
 \ar[u]^-{(\chi,\phi)^*}
 \ar[d]_-{\re^*} &
 \scriptstyle{[\^e\^f-\delta\^g]}
 \ar@{|->}[u]
 \ar@{|->}[d] \\
 H^{m+n}(S^m\cro S^n).
 &
 \scriptstyle{[\^e][\^f]}
 }
 $$
 Collating it with \eqref{homol} and
 using \eqref{int},
 we get
 $$
 \([q],[S^{m+n-1}]\)=
 (-1)^{m+n}
 \([\^e][\^f],[S^m\cro S^n]\)=
 (-1)^{mn+m+n}\Delta.
 $$
 This is what we need because
 $h(\[\+u,\+v\])=h([a])=[q]$.
 \qed
 \end {demo}

 Let $\Gamma$ be an open cover of $X$.
 Consider the differential graded ring $C^\bullet(\Gamma)$
 of $\Gamma$-cochains of $X$
 (that is,
 functions on the set of singular simplices
 subordinate to $\Gamma$).
 The projection
 $$
 ?|_\Gamma:C^\bullet(X)\to C^\bullet(\Gamma)
 $$
 is a morphism of differential graded rings;
 it induces an isomorphism of cohomology rings,
 $$
 ?|_\Gamma:H^\bullet(X)\to H^\bullet(\Gamma).
 $$
 
 \begin {claim} [9.3. Lemma.]
 Given $a\in Y^X_{e,f}$,
 choose $\~p\in C^{m-1}(\Gamma)$ such that
 $\delta\~p=a^\#(e)|_\Gamma$
 and
 put
 $$
 \~q=
 \~pa^\#(f)|_\Gamma-a^\#(g)|_\Gamma\in
 C^{m+n-1}(\Gamma).
 $$
 Then
 $\delta\~q=0$
 and
 $h([a])|_\Gamma=[\~q]$ in $H^{m+n-1}(\Gamma)$.
 \qed
 \end {claim}

 We suppose that
 $X$ and $Y$ are cellular spaces
 and
 $X$ is compact.

 \begin {claim} [9.4. Theorem.]
 Consider an ensemble $A\in\<Y^X\>$,
 $$
 A=
 \sum_i
 u_i\`a_i\',
 $$
 where $a_i\in Y^X_{e,f}$,
 such that
 $A\%=20$.
 Then
 $$
 \sum_i
 u_ih([a_i])
 =0
 $$
 in $H^{m+n-1}(X)$.
 \end {claim}

 Thus
 $h$ may be called a {\it partial\/} invariant of order
 at most $2$.

 \begin {demo} [Proof.]
 Using Corollary~6.2,
 replace $a_i$ by homotopic maps so that
 $A\$=2\Gamma0$
 for some open cover $\Gamma$ of $X$.

 Let $B\subseteq C^m(\Gamma)$ be the subgroup generated by
 the coboundaries $a_i^\#(e)|_\Gamma$.
 It is free
 because
 finitely generated
 and 
 torsion-free.
 Thus
 there is a homomorphiam $P:B\to C^{m-1}(\Gamma)$ such that
 $\delta P(b)=b$, $b\in B$.
 Put
 $$
 \~q_i=
 P(a^\#(e)|_\Gamma)a^\#(f)|_\Gamma-a^\#(g)|_\Gamma\in
 C^{m+n-1}(\Gamma).
 $$
 By Lemma~9.3,
 $\delta\~q_i=0$
 and
 $$
 h([a_i])|_\Gamma=[\~q_i]
 $$
 in $H^{m+n-1}(\Gamma)$.

 Take a singular simplex $\sigma:\Delta^{m+n-1}\to G$,
 $G\in\Gamma$.
 Let
 $$
 \sigma':\Delta^{m-1}\to G
 \quad
 \text{and}
 \quad
 \sigma'':\Delta^n\to G
 $$
 be its front and back faces,
 respectively.

 The group $\Hom(B,\Q)$ is formed by homomorphisms $\(?,T\)$,
 where $T$ runs over $C_m(\Gamma;\Q)$,
 the group of rational $\Gamma$-chains in $X$.
 Thus
 there is a chain $T\in C_m(\Gamma;\Q)$ such that
 $$
 \(P(b),\sigma'\)=
 \(b,T\),
 \qquad
 b\in B.
 $$
 We have
 $$
 T=
 \sum_k
 c_k\tau_k,
 $$
 where
 $c_k\in\Q$
 and
 $\tau_k:\Delta^m\to G_k$,
 $G_k\in\Gamma$.
 Thus
 $$
 \(P(a_i^\#(e)|_\Gamma),\sigma'\)=
 \(a_i^\#(e)|_\Gamma,T\)=
 \sum_k
 c_k\(a_i^\#(e)|_\Gamma,\tau_k\).
 $$
 We get
 \begin {multline*}
 \(\~q_i,\sigma\)=
 (-1)^{(m-1)n}
 \(P(a_i^\#(e)|_\Gamma),\sigma'\)
 \(a_i^\#(f)|_\Gamma,\sigma''\)-
 \(a_i^\#(g)|_\Gamma,\sigma\)= \\ =
 (-1)^{(m-1)n}
 \sum_k
 c_k
 \(a_i^\#(e)|_\Gamma,\tau_k\)
 \(a_i^\#(f)|_\Gamma,\sigma''\)-
 \(a_i^\#(g)|_\Gamma,\sigma\)= \\ =
 (-1)^{(m-1)n}
 \sum_k
 c_k
 \((a_i|_{G\cup G_k})^\#(e),\tau_k\)
 \((a_i|_{G\cup G_k})^\#(f),\sigma''\)-
 \((a_i|_G)^\#(g),\sigma\).
 \end {multline*}
 We have found functions
 $R_k:Y^{G\cup G_k}\to\Q$
 and
 $S:Y^G\to\Q$
 such that
 $$
 \(\~q_i,\sigma\)=
 \sum_k
 R_k(a_i|_{G\cup G_k})-
 S(a_i|_G)
 $$
 for all $i$.
 Since $A\$=2\Gamma0$,
 we have
 $A|_{G\cup G_k}=0$
 and
 $A|_G=0$.
 Thus
 $$
 \sum_i
 u_i\(\~q_i,\sigma\)=
 0.
 $$

 Since $\sigma$ was taken arbitrarily,
 we have
 $$
 \sum_i
 u_i\~q_i=
 0.
 $$
 We get
 $$
 \sum_i
 u_ih([a_i])|_\Gamma=
 \sum_i
 u_i[\~q_i]=
 0.
 $$
 Since restriction to $\Gamma$ here is an isomorphism,
 we get
 $$
 \sum_i
 u_ih([a_i])=
 0.
 \QED
 $$
 \end {demo}

 \begin {claim} [9.5. Corollary.]
 Let $a,b\in Y^X_{e,f}$ satisfy $a\%\sim2b$.
 Then
 $h([a])=h([b])$.
 \end {claim}

 \begin {demo} [Proof.]
 There is an ensemble $A\in\<Y^X\>$,
 $$
 A=
 \sum_i
 u_i\`a_i\',
 $$
 where $a_i\sim a$,
 such that
 $A\%=2\`b\'$.
 Since $A\|={\{\0\}}\`b\'$,
 we have
 $$
 \sum_i u_i=1.
 $$
 By Theorem~9.4,
 $$
 \sum_i
 u_ih([a_i])=
 h([b]).
 $$
 Since $[a_i]=[a]$,
 we get $h([a])=h([b])$.
 \qed
 \end {demo}


 \head {\S~10. Maps of $S^p\cro S^n$}


 This section does not depend of the rest of the paper.
 We
 recall a theorem of G.~W.~Whitehead
 about the fibration of free spheroids
 (Theorem~10.1)
 and
 deduce Lemma~10.3 about certain maps $S^p\cro S^n\to Y$
 (we need it in \S~11).

 We fix
 numbers $p,n\ge1$
 and
 a space $Y$.
 Let $\Omega^nY$ be the space of maps $S^n\to Y$,
 as usual.
 Let
 $$
 \epsilon:
 S^p\cro S^n\to
 S^p\sma S^n\to
 S^{p+n}
 $$
 be the composition of
 the projection and
 the standard homeomorphism.
 For a map $w:S^{p+n}\to Y$,
 introduce the map
 $$
 \nabla^n(w):S^p\to\Omega^nY,
 \qquad
 \nabla^n(w)(t)(z)=w(\epsilon(t,z)).
 $$
 Introduce the isomorphism
 $$
 \+\nabla^n:\pi_{p+n}(Y)\to\pi_p(\Omega^nY),
 \qquad
 [w]\mapsto[\nabla^n(w)].
 $$

 Let
 $$
 \mu:S^n\to S^n\bou S^n
 $$
 be the standard comultiplication.
 Consider the usual multiplication
 $$
 \Omega^nY\cro\Omega^nY
 \xto\#
 \Omega^nY,
 \qquad
 v_1\#v_2:
 S^n
 \xto\mu
 S^n\bou S^n
 \xto{v_1\Bou v_2}
 Y.
 $$
 For a map $v:S^n\to Y$,
 introduce the map
 $$
 \tau_v:
 \Omega^nY
 \xto{v\#?}
 (\Omega^nY,v\#\0),
 $$
 where the target is $\Omega^nY$ with the specified new basepoint.
 It induces the isomorphism
 $$
 \tau_{v\,*}:
 \pi_p(\Omega^nY)
 \to
 \pi_p(\Omega^nY,v\#\0).
 $$

 Let $\Lambda^nY$ be the space of {\it unbased\/} maps $S^n\to Y$.
 Consider the fibration
 $$
 \rho:\Lambda^nY\to Y,
 \qquad
 v\mapsto v(\0).
 $$
 We have $\rho^{-1}(\0)=\Omega^n(Y)$.
 
 \begin {claim} [10.1. Theorem \rm (G.~W.~Whitehead).]
 For a map $v:S^n\to Y$,
 the composition
 $$
 \+\Gamma:
 \pi_{p+1}(Y)
 \xto{\[?,[v]\]}
 \pi_{p+n}(Y)
 \xto{\+\nabla^n}
 \pi_p(\Omega^nY)
 \xto{\tau_{v\,*}}
 \pi_p(\Omega^nY,v\#\0)
 $$
 coincides up to a sign with
 the connecting homomorphism
 of the fibration $\rho$
 at the point $v\#\0\in\Omega^nY$.
 Consequently,
 the composition
 $$
 \pi_{p+1}(Y)
 \xto{\+\Gamma}
 \pi_p(\Omega^nY,v\#\0)
 \xto{\ins_*}
 \pi_p(\Lambda^nY,v\#\0)
 $$
 is zero.
 \end {claim}

 \begin {dem}
 See
 \cite[Theorem~(3.2)]{WhiG}
 and
 \cite[\S~3]{WhiH}.
 \qed
 \end {dem}

 For a map $v:S^n\to Y$,
 introduce the homomorphism
 $$
 \+\Psi_v:
 \pi_{p+n}(Y)
 \xto{\+\nabla^n}
 \pi_p(\Omega^nY)
 \xto{\tau_{v\,*}}
 \pi_p(\Omega^nY,v\#\0)
 \xto{\ins_*}
 \pi_p(\Lambda^nY,v\#\0).
 $$
 By Theorem~10.1,
 \begin {equation} \label {Psi}
 \+\Psi_v(\[\+u,[v]\])=0,
 \qquad
 \+u\in\pi_{p+1}(Y).
 \end {equation}
 For maps
 $v:S^n\to Y$
 and
 $w:S^{p+n}\to Y$,
 introduce the map
 $$
 \Psi_v(w):
 S^p
 \xto{\nabla^n(w)}
 \Omega^nY
 \xto{\tau_v}
 (\Omega^nY,v\#\0)
 \xto\inc
 (\Lambda^nY,v\#\0).
 $$
 Clearly,
 $$
 [\Psi_v(w)]=\+\Psi_v([w])
 $$
 in $\pi_p(\Lambda^nY,v\#\0)$.

 Introduce the map
 \begin {equation}
 \Phi:
 S^p\cro S^n
 \xto{\id\cro\mu}
 S^p\cro(S^n\bou S^n)
 \xto\theta
 S^n\bou S^{p+n},
 \end {equation}
 where
 $$
 \theta:
 \quad
 (t,\ins_1(z))\mapsto\ins_1(z),
 \quad
 (t,\ins_2(z))\mapsto\ins_2(\epsilon(t,z)),
 \qquad
 t\in S^p,
 \quad
 z\in S^n.
 $$
 For maps
 $v:S^n\to Y$
 and
 $w:S^{p+n}\to Y$,
 introduce the map
 \begin {equation} \label {Xi}
 \Xi(v,w):
 S^p\cro S^n
 \xto\Phi
 S^n\bou S^{p+n}
 \xto{v\Bou w}
 Y.
 \end {equation}
 For elements
 $\+v\in\pi_n(Y)$
 and
 $\+w\in\pi_{p+n}(Y)$,
 put
 \begin {equation} \label {XXi}
 \+\Xi(\+v,\+w)=[\Xi(v,w)]
 \in
 [S^p\cro S^n,Y],
 \end {equation}
 where $v$ and $w$ are representatives of $\+v$ and $\+w$,
 respectively.

 For maps
 $v_0:S^n\to Y$
 and
 $V:S^p\to(\Lambda^nY,v_0)$,
 introduce the map
 $$
 V^\c:S^p\cro S^n\to Y,
 \qquad
 (t,z)\mapsto V(t)(z).
 $$
 For $\+V\in\pi_p(\Lambda^nY,v_0)$,
 put
 $$
 \+V^{\+\c}=[V^\c]\in[S^p\cro S^n,Y],
 $$
 where $V$ is a representative of $\+V$.

 \begin {claim} [10.2. Lemma.]
 For maps
 $v:S^n\to Y$
 and
 $w:S^{p+n}\to Y$,
 one has
 $$
 \Xi(v,w)=\Psi_v(w)^\c:
 S^p\cro S^n\to Y.
 $$
 Consequently,
 $$
 \+\Xi([v],[w])=\+\Psi_v([w])^{\+\c}
 $$
 in $[S^p\cro S^n,Y]$.
 \end {claim}

 \begin {demo} [Proof.]
 Take a point $(t,z)\in S^p\cro S^n$.
 We have $\mu(z)=\ins_k(\~z)$
 in $S^n\bou S^n$
 for some
 $k\in\{1,2\}$ and
 $\~z\in S^n$.
 We have
 \begin {align*}
 \theta(t,\mu(z))=
 \theta(t,\ins_k(\~z))=
 \qquad
 \text{(if $k=1$)}
 \qquad
 &=
 \ins_1(\~z), \\
 \text{(if $k=2$)}
 \qquad
 &=
 \ins_2(\epsilon(t,\~z))
 \end {align*}
 in $S^n\bou S^{p+n}$.
 Thus
 \begin {multline*}
 \Xi(v,w)(t,z)=
 ((v\Bou w)\circ\Phi)(t,z)= \\
 =
 ((v\Bou w)\circ\theta\circ(\id\cro\mu))(t,z)=
 (v\Bou w)(\theta(t,\mu(z)))= \\
 \begin {aligned}
 \text{(if $k=1$)}
 \qquad
 &=
 (v\Bou w)(\ins_1(\~z))=
 v(\~z), \\
 \text{(if $k=2$)}
 \qquad
 &=
 (v\Bou w)(\ins_2(\epsilon(t,\~z)))=
 w(\epsilon(t,\~z)).
 \end {aligned}
 \end {multline*}
 On the other hand,
 \begin {multline*}
 \Psi_v(w)^\c(t,z)=
 \Psi_v(w)(t)(z)=
 \tau_v(\nabla^n(w)(t))(z)= \\
 =
 (v\#\nabla^n(w)(t))(z)=
 (v\Bou\nabla^n(w)(t))(\mu(z))=
 (v\Bou\nabla^n(w)(t))(\ins_k(\~z))= \\
 \begin {aligned}
 \text{(if $k=1$)}
 \qquad
 &=
 v(\~z), \\
 \text{(if $k=2$)}
 \qquad
 &=
 \nabla^n(w)(t)(\~z)=
 w(\epsilon(t,\~z)).
 \end {aligned}
 \end {multline*}
 The same.
 \qed
 \end {demo}

 \begin {claim} [10.3. Lemma.]
 For elements
 $\+u\in\pi_{p+1}(Y)$,
 $\+v\in\pi_n(Y)$,
 and
 $\+w\in\pi_{p+n}(Y)$,
 one has
 $$
 \+\Xi(\+v,\[\+u,\+v\]+\+w)=\+\Xi(\+v,\+w)
 $$
 in $[S^p\cro S^n,Y]$.
 \end {claim}

 \begin {demo} [Proof.]
 Choose a representative $v:S^n\to Y$ of $\+v$.
 By \eqref{Psi},
 $$
 \+\Psi_v(\[\+u,\+v\]+\+w)=\+\Psi_v(\+w)
 $$
 in $\pi_p(\Lambda^nY,v\#\0)$.
 Applying Lemma~10.2 yields the desired equality.
 \qed
 \end {demo}

 For a map
 $w:S^{p+n}\to Y$,
 introduce the map
 $$
 \xi(w):
 S^p\cro S^n
 \xto\epsilon
 S^{p+n}
 \xto{w}
 Y.
 $$
 For an element
 $\+w\in\pi_{p+n}(Y)$,
 put
 \begin {equation} \label {xxi}
 \+\xi(\+w)=[\xi(w)]
 \in
 [S^p\cro S^n,Y],
 \end {equation}
 where $w$ is a representative of $\+w$.

 \begin {claim} [10.4. Lemma.]
 For en element $\+w\in\pi_{p+n}(Y)$,
 one has
 $$
 \+\Xi(0,\+w)=\+\xi(\+w)
 $$
 in $[S^p\cro S^n,Y]$.
 \end {claim}

 \begin {demo} [Proof.]
 Choose a representative $w:S^{p+n}\to Y$ of $\+w$.
 Consider the diagram
 $$
 \xymatrix {
 S^p\cro(S^n\bou S^n)
 \ar[rrrr]^-{\theta}
 \ar[dddd]_-{\id\cro(\0\Bou\id)}
 &&&&
 S^n\bou S^{p+n}
 \ar[ddl]^>>>>>>{\0\Bou w}
 \ar[dddd]^-{\0\Bou\id} \\ \\
 &
 S^p\cro S^n
 \ar[uul]_-{\id\cro\mu}
 \ar[uurrr]^-{\Phi}
 \ar@/^1ex/[rr]^>>>>{\Xi(\0,w)}
 \ar@/_1ex/[rr]_-{\xi(w)}
 \ar[ddl]^-{\id}
 \ar[ddrrr]_-{\epsilon}
 &&
 Y
 & \\ \\
 S^p\cro S^n
 \ar[rrrr]_-{\epsilon}
 &&&&
 S^{p+n}.
 \ar[uul]_-{w}
 }
 $$
 Since the map
 $$
 S^n
 \xto\mu
 S^n\bou S^n
 \xto{\0\Bou\id}
 S^n
 $$
 is homotopic to the identity,
 the left triangle is homotopy commutative.
 The other empty triangles and the square are commutative.
 It follows that the parallel curved arrows are homotopic.
 \qed
 \end {demo}


 \head {\S~11. Fineness of $2$-similarity}


 Put $X=S^p\cro S^n$
 ($p\ge1$, $n\ge2$).
 Let $Y$ be a space with
 elements
 $\+u\in\pi_{p+1}(Y)$ and
 $\+v\in\pi_n(Y)$.
 Consider
 the Whitehead product $\[\+u,\+v\]\in\pi_{p+n}(Y)$ and
 the homotopy classes
 $$
 \+k(t)=\+\xi(t\[\+u,\+v\])\in[X,Y],
 \qquad
 t\in\Z
 $$
 (see \eqref{xxi}).

 \begin {claim} [11.1. Lemma.]
 Let
 $L$ be an abelian group
 and
 $f:[X,Y]\to L$ be an invariant of order at most $r$.
 Then
 $$
 f(\+k(r!+t))=f(\+k(t)),
 \qquad
 t\in\Z.
 $$
 \end {claim}

 \begin {demo} [Proof \rm (after {\cite[Lemma~1.5]{uns}}).]
 We will use the homotopy classes
 $$
 \+K(s,t)=\+\Xi(s\+v,t\[\+u,\+v\])\in[X,Y],
 \qquad
 s,t\in\Z
 $$
 (see \eqref{XXi}).
 By Lemma~10.4,
 \begin {equation} \label {Xixi}
 \+K(0,t)=\+k(t).
 \end {equation}
 We have
 \begin {equation} \label {period}
 \+K(s,m+t)=\+K(s,t)
 \qquad
 \text {if $s\divs m$}
 \end {equation}
 because
 \begin {multline*}
 \+\Xi(s\+v,(m+t)\[\+u,\+v\])=
 \+\Xi(s\+v,\[(m/s)\+u,s\+v\]+t\[\+u,\+v\])= \\
 \text{(by Lemma~10.3)}
 \qquad
 =
 \+\Xi(s\+v,t\[\+u,\+v\]).
 \end {multline*}

 Consider the wedge of
 $r$ copies of $S^n$
 and
 two copies of $S^{p+n}$
 $$
 W=S^n\bou\dotso\bou S^n\bou S^{p+n}\bou S^{p+n}
 $$
 and the maps
 $$
 \Lambda(d)=
 \lambda_1(d_1)\bou\dotso\lambda_r(d_r)\bou
 \lambda_{r+1}(d_{r+1})\bou\id:
 W\to W,
 $$
 $d=(d_1,\dotsc,d_{r+1})\in\E^{r+1}$,
 as in \S~3.
 Put
 $$
 \mu=\mu_1\bou\mu_2:
 S^n\bou S^{p+n}
 \to W,
 $$
 where
 $$
 \mu_1:S^n\to S^n\bou\dotso\bou S^n
 \quad
 \text{and}
 \quad
 \mu_2:S^{p+n}\to S^{p+n}\bou S^{p+n}
 $$
 are the comultiplications.
 Choose a map $q:W\to Y$ with
 $$
 [q]=\+v\Bou\dotso\Bou\+v\Bou r!\[\+u,\+v\]\Bou t\[\+u,\+v\].
 $$
 Consider the ensemble $A\in\<Y^X\>$,
 $$
 A=\sum_{d\in\E^{r+1}}(-1)^{|d|}\`a(d)\',
 $$
 where
 $$
 a(d):
 X
 \xto{\Phi}
 S^n\bou S^{p+n}
 \xto{\mu}
 W
 \xto{\Lambda(d)}
 W
 \xto{q}
 Y,
 $$
 where $\Phi$ is as in \eqref{Xi}.
 By Lemma~3.1,
 $A\%=r0$.
 Clearly,
 $$
 [q\circ\Lambda(d)\circ\mu]=
 (d_1+\dotso d_r)\+v\Bou
 (d_{r+1}r!+t)\[\+u,\+v\]
 $$
 in $[S^n\bou S^{p+n},Y]$.
 Thus,
 by the construction of $\+K(s,t)$,
 $$
 [a(d)]=
 \+K(d_1+\dotso d_r,d_{r+1}r!+t)
 $$
 in $[X,Y]$.
 Thus,
 since $f$ has order at most $r$,
 $$
 \sum_{d\in\E^{r+1}}
 (-1)^{|d|}f(\+K(d_1+\dotso d_r,d_{r+1}r!+t))=0.
 $$
 By \eqref{period},
 $\+K(d_1+\dotso d_r,d_{r+1}r!+t)$
 does not depend on $d_{r+1}$
 if $(d_1,\dotsc,d_r)\neq(0,\dotsc,0)$.
 Thus
 the corresponding summands cancel out.
 We get $f(\+K(0,t))-f(\+K(0,r!+t))=0$.
 By \eqref{Xixi},
 this is what we need.
 \qed
 \end {demo}

 Let classes 
 $E\in H^{p+1}(Y)$ and
 $F\in H^n(Y)$
 satisfy
 $EF=0$ in $H^{p+n+1}(Y)$.
 Put,
 as in Lemma~9.2,
 $$
 \Delta=
 \(\+u^*(E),[S^{p+1}]\)
 \(\+v^*(F),[S^n]\)+
 (-1)^{(p+1)n}
 \(\+u^*(F),[S^{p+1}]\)
 \(\+v^*(E),[S^n]\)
 \in\Z.
 $$

 If $Y=S^{p+1}\bou S^n$ with
 $\+u=[\ins_1]$ and
 $\+v=[\ins_2]$,
 taking obvious $E$ and $F$ yields $\Delta=1$.
 If
 $p=n-1$ and
 $Y=S^n$ with
 $\+u=\+v=[\id]$,
 taking obvious equal $E$ and $F$ yields $\Delta=1+(-1)^n$.

 \begin {claim} [11.2. Lemma.]
 If $\Delta\ne0$,
 the classes $\+k(t)$, $t\in\Z$, are pairwise not $2$-similar.
 \end {claim}

 \begin {demo} [Proof.]
 Choose
 cocycles
 $e\in C^{p+1}(Y)$ and
 $f\in C^n(Y)$
 representing
 $E$ and
 $F$,
 respectively.
 Choose a cochain $g\in C^{p+n}(Y)$ with $\delta g=ef$.
 Consider the corresponding Hopf invariants
 (see \S~9)
 $$
 h_0:\pi_{p+n}(Y)\to
 H^{p+n}(S^{p+n})
 \quad
 \text{and}
 \quad
 h:[X,Y]_{e,f}\to
 H^{p+n}(X).
 $$
 By Lemma~9.2,
 $$
 \(h_0(\[\+u,\+v\]),[S^{p+n}]\)=
 (-1)^{pn+p+1}\Delta.
 $$
 We have the decomposition
 $$
 \xymatrix {
 \+k(t):
 X
 \ar[r]^-{\epsilon} &
 S^{p+n}
 \ar@{~>}[r]^-{t[\id]} &
 S^{p+n}
 \ar@{~>}[r]^-{\[\+u,\+v\]} &
 Y
 }
 $$
 (the wavy arrows denote homotopy classes).
 Clearly,
 $\+k(t)\in[X,Y]_{e,f}$.
 Since
 the Brouwer degree of $\epsilon$ is $1$ and
 that of $t[\id]$ is $t$,
 Lemma~9.1 yields
 $$
 \(h(\+k(t)),[X]\)=
 (-1)^{pn+p+1}\Delta t.
 $$
 By Corollary~9.5,
 the classes $\+k(t)$, $t\in\Z$, are pairwise
 not $2$-similar
 if $\Delta\ne0$.
 \qed
 \end {demo}

 \subhead {Moral.}
 Suppose that
 $\Delta\ne0$.
 The classes
 $\+k(0)$ ($=[\0]$) and
 $\+k(2)$
 in $[X,Y]$,
 which are not $2$-similar
 by Lemma~11.2,
 cannot be distinguished by an invariant of order at most $2$
 by Lemma~11.1.
 Recall that
 $(X,Y)$ can be
 $(S^p\cro S^n,S^{p+1}\bou S^n)$
 for any $p\ge1$ and $n\ge2$ or
 $(S^{n-1}\cro S^n,S^n)$
 for even $n\ge2$.



 {\noindent \tt ssp@pdmi.ras.ru}

 {\noindent \tt http://www.pdmi.ras.ru/\til{}ssp}

 \end {document}